\documentclass{amsart}

\usepackage{mathrsfs}
\usepackage{mathrsfs}
\usepackage{amsfonts}
\usepackage{amsfonts}
\usepackage{mathrsfs}
\usepackage{amsfonts}
\usepackage{amsfonts}
\usepackage{amsfonts}
\usepackage{mathrsfs}
\usepackage{amsfonts}
\usepackage{amssymb,amsmath}
\usepackage{amsthm}
\usepackage{amsfonts}

\newtheorem{theorem}{Theorem}[section]
\newtheorem{lemma}[theorem]{Lemma}

\theoremstyle{definition}

\theoremstyle{remark}

\numberwithin{equation}{section}

\begin{document}

\title[]
 {H$\ddot{o}$lder continuity for stochastic fractional heat equation with colored noise}

\author[]{Kexue Li}
\address{School of Mathematics and Statistics,
Xi'an Jiaotong University,
 Xi'an
710049, China}

\email{kexueli@gmail.com}





\begin{abstract}
In this paper, we consider semilinear stochastic fractional heat equation $\frac{\partial}{\partial t}u_{\beta,t}(x)=\triangle^{\alpha/2}u_{\beta,t}(x)+\sigma(u_{\beta,t}(x))\eta_{\beta}$. The Gaussian noise $\eta_{\beta}$ is assumed to be colored in space with covariance of the form $E(\eta_{\beta}(t,x)\eta_{\beta}(s,y))=\delta(t-s)f_{\beta}(x-y)$, where $f_{\beta}$ is the Riesz kernel $f_{\beta}(x)\propto |x|^{-\beta}$. We obtain the spatial and temporal H$\ddot{\mbox{o}}$lder continuity of the mild solution.

\end{abstract}

\maketitle


\section{Introduction}
In this paper, we consider the following stochastic fractional heat equation
\begin{equation}\label{sfh}
\ \left\{\begin{aligned}
&\frac{\partial}{\partial t}u_{\beta, t}(x)=\triangle^{\alpha/2}u_{\beta, t}(x)+\sigma(u_{\beta,t}(x))\eta_{\beta}
&,\ t>0,\ x\in R,\\
&u_{\beta,0}(x)=\phi(x),
\end{aligned}\right.
\end{equation}
where $1<\alpha\leq2$,  $\triangle^{\alpha/2}:=-(-\Delta)^{\alpha/2}$ denotes the fractional Laplacian defined by the Fourier transform
\begin{align*}
(\mathcal{F}(-\triangle)^{\alpha/2}u)(\xi)=(2\pi|\xi|)^{\alpha}\mathcal{F}(u)(\xi),
\end{align*}
here $\mathcal{F}$ denotes the Fourier transform,
\begin{align}\label{Fourier transform}
(\mathcal{F}\varphi)(\xi)=\int_{\mathbb{R}}e^{-2i\pi\xi x}\varphi(x)dx.
\end{align}
$\eta_{\beta}$ is the Gaussian space time colored noise with covariance of the form
\begin{equation}
E[\eta_{\beta}(t,x)\eta_{\beta}(s,y)]=\delta(t-s)f_{\beta}(x-y),
\end{equation}
 where (\cite{Robert}, Ex.1)
\begin{equation}\label{transform}
f_{\beta}(x)=c_{1-\beta}g_{\beta}(x)=(\mathcal{F}g_{1-\beta})(x), \ g_{\beta}(x)=\frac{1}{|x|^{\beta}}, \ \beta\in (0,1),
\end{equation}
and
\begin{equation}
c_{\beta}=\frac{2\sin(\beta\pi/2)\Gamma(1-\beta)}{(2\pi)^{1-\beta}},
\end{equation}
where $\Gamma(\cdot)$ is the Gamma function. \\
We assume that the following conditions hold: \\
(A1) $\phi$ is bounded and $\rho$-H$\ddot{\mbox{o}}$lder continuous. \\
(A2) $\sigma$ is Lipschitz continuous and there exists a constant $K$ such that $|\sigma(x)-\sigma(y)|\leq K|x-y|$ and $|\sigma(x)|\leq K(1+|x|)$.

The mild solutions are the solutions of the integral equations
\begin{align}\label{mild solution}
u_{\beta,t}(y)=(u_{\beta,0}\ast p_{t})(y)+\int_{0}^{t}\int_{\mathbb{R}}p_{t-s}(x-y)\sigma(u_{\beta,s}(x))\eta_{\beta}(ds,dx).
\end{align}
where the fractional heat kernel $p_{t}(x)$ is the fundmental solution of
\begin{align}\label{fundmental}
v_{t}=\Delta^{\alpha/2}v,
\end{align}
and $\ast$ denotes the usual convolution operator,  $(f\ast g)(x)=\int_{\mathbb{R}}f(x-y)g(y)dy$.
Since  $0<\beta<1$, we can get the existence and uniqueness of the mild solution of (\ref{sfh}) (see, e.g., \cite{Foondun,MM}).
It is known that $p_{t}(x)$ satisfies the following inequality (\cite{TG, CKS, KT})
\begin{align}\label{alpha}
\frac{c_{1}t}{(t^{1/\alpha}+|x|)^{1+\alpha}}\leq p_{t}(x)\leq \frac{c_{2}t}{(t^{1/\alpha}+|x|)^{1+\alpha}},
\end{align}
where $t>0$, $x\in \mathbb{R}$, $c_{1}$ and $c_{2}$ are positive constants depending on $\alpha$.

In the very recent paper \cite{PB}, Bezdek considered the following equations
\begin{equation}\label{colored}
\ \left\{\begin{aligned}
&\frac{\partial}{\partial t}u_{\beta, t}(x)=\frac{\kappa}{2}\triangle u_{\beta, t}(x)+\sigma(u_{\beta,t}(x))\eta_{\beta}
&,\ t> 0,\ x\in R,\\
&u_{\beta,0}(x)=\phi(x),
\end{aligned}\right.
\end{equation}
where $\kappa>0$ and $\eta_{\beta}$ is the Gaussian noise colored in space and white in time. Stochastic PDEs with colored noise has been studied in many papers (see, e.g., \cite{BC,BA,Conus,Robert}).

Bezdek has obtained the H$\ddot{\mbox{o}}$lder continuity estimates which take into account $\beta$ as a variable, the results are novel in that sense.  In this paper, based on some estimates of the fractional heat kernel, we will show the spatial and temporal H$\ddot{\mbox{o}}$lder continuity for the mild solution of stochastic fractional heat equations (\ref{sfh}).


\section{The spatial and temporal H$\ddot{\mbox{o}}$lder continuity}
\subsection{Some lemma} In this subsection, we will prove some lemmas, which will be used in next subsections.
 We use $C$ to denote generic constant, which may change from line to line.

\begin{lemma}\label{space}
For all $t>0$ and $x\in \mathbb{R}$,
\begin{align}
\int_{\mathbb{R}}|p_{t}(y-x)-p_{t}(y)|dy\leq C(\frac{|x|}{t^{1/\alpha}}\wedge 1),
\end{align}
where $C$ does not depend on $t$ or $x$.
\end{lemma}

\begin{proof}
For all $r>0$,
define
\begin{align}
\mu(r)=\mu(r,t):=\sup_{z\in \mathbb{R}, |z|\leq r}\int_{\mathbb{R}}|p_{t}(y-z)-p_{t}(y)|dy.
\end{align}
Then
\begin{align}\label{derivative}
\mu(|x|)=\sup_{z\in (0,|x|)}\int_{-\infty}^{\infty}\big|\int_{y-z}^{y}\frac{\partial p_{t}(\xi)}{\partial \xi}d\xi\big|dy.
\end{align}
By (2.3) of \cite{TG} (or Lemma 5 in \cite{KT}), we have
\begin{align}\label{gradient}
\big|\frac{\partial p_{t}(\xi)}{\partial \xi}\big|\leq C\frac{t|\xi|}{(t^{1/\alpha}+|\xi|)^{3+\alpha}},
\end{align}
where $C$ only depends on $\alpha$.\\
Taking (\ref{gradient}) into (\ref{derivative}) to get
\begin{align}
\mu(|x|)&\leq C|x|\int_{-\infty}^{\infty}\frac{t|w|}{(t^{1/\alpha}+|\xi|)^{3+\alpha}}d\xi\nonumber\\
&=\frac{C|x|}{t^{1/\alpha}}\int_{-\infty}^{\infty}\frac{|\nu|}{(1+|\nu|)^{3+\alpha}}d\nu\nonumber\\
&\leq \frac{C|x|}{t^{1/\alpha}}.
\end{align}
On the other hand, since $|p_{t}(y-x)-p_{t}(y)|\leq p_{t}(y-x)+p_{t}(y)$ and $\int_{\mathbb{R}}p_{t}(y)dy=1$, we have $\mu(|x|)\leq 2$.
\end{proof}
\begin{lemma}\label{time}
For all $t, \varepsilon>0$, we have
\begin{align}
\int_{\mathbb{R}}|p_{t+\varepsilon}(y)-p_{t}(y)|dy\leq C(\log(t+\varepsilon)-\log(t))\wedge 1).
\end{align}
\end{lemma}

\begin{proof}
From (\ref{fundmental}) and Proposition 2.1 in \cite{Vazquez}, it is easy to show that
\begin{align}
|\frac{\partial p_{t}(y)}{\partial t}|\leq \frac{Cp_{t}(y)}{t}.
\end{align}
Then we have
\begin{align}
&\int_{\mathbb{R}}|p_{t+\varepsilon}(y)-p_{t}(y)|dy=\int_{\mathbb{R}}\big|\int_{t}^{t+\varepsilon}\frac{\partial p_{s}(y)}{\partial s}ds\big|dy\leq C\int_{\mathbb{R}}\int_{t}^{t+\varepsilon}\frac{p_{s}(y)}{s}dsdy\nonumber\\
&=C\int_{t}^{t+\epsilon}\frac{1}{s}ds=C(\log(t+\varepsilon)-\log(t)).
\end{align}
On the other hand, we have $\int_{\mathbb{R}}|p_{t+\varepsilon}(y)-p_{t}(y)|dy\leq 2$.
\end{proof}

\begin{lemma}\label{holder}
Let $0<\rho<1$ and let $w$ be a bounded $\rho$-H$\ddot{\mbox{o}}$lder continuous function, then there exists $C>0$ such that for every $t>0$, $\delta>0$, $x\in \mathbb{R}$, $z\in \mathbb{R}$, we have
\begin{align*}
\big|\int_{\mathbb{R}}(p_{t}(x-y)-p_{t}(z-y))w(y)dy\big|&\leq C|x-z|^{\rho}, \\
\big|\int_{\mathbb{R}}(p_{t+\delta}(x-y)-p_{t}(x-y))w(y)dy\big|&\leq C\delta^{\rho/\alpha}.
\end{align*}
\end{lemma}

\begin{proof}
\begin{align*}
&\big|\int_{\mathbb{R}}(p_{t}(x-y)-p_{t}(z-y))w(y)dy\big|=\big|\int_{\mathbb{R}}p_{t}(z-y)w(y+x-z)dy-\int_{\mathbb{R}}p_{t}(z-y)w(y)dy\big|\\
&=\big|\int_{\mathbb{R}}p_{t}(z-y)(w(y+x-z)-w(y))dy\big|\leq C|x-z|^{\rho}\int_{\mathbb{R}}p_{t}(z-y)dy=C|x-z|^{\rho}.
\end{align*}
By the semigroup property of $p_{t}$ and note that $\int_{\mathbb{R}}p_{\delta}(y)dy=\int_{\mathbb{R}}p_{t}(x-z)dz=1$,
\begin{align}\label{delta}
&\big|\int_{\mathbb{R}}(p_{t+\delta}(x-y)-p_{t}(x-y))w(y)dy\big|=\big|\int_{\mathbb{R}}\big(\int_{\mathbb{R}}p_{t}(x-z)p_{\delta}(z-y)dz\big)w(y)dy-\int_{\mathbb{R}}p_{t}(x-y)w(y)dy\big|\nonumber\\
&=\big|\int_{\mathbb{R}}\big(\int_{\mathbb{R}}p_{t}(x-z)p_{\delta}(y)dz\big)w(z-y)dy-\int_{\mathbb{R}}p_{t}(x-z)w(z)dz\big|\nonumber\\
&=\big|\int_{\mathbb{R}}p_{\delta}(y)\big(\int_{\mathbb{R}}p_{t}(x-z)w(z-y)dz\big)dy-\int_{\mathbb{R}}p_{\delta}(y)\big(\int_{\mathbb{R}}p_{t}(x-z)w(z)dz\big)dy\big|\nonumber\\
&\leq \int_{\mathbb{R}}p_{\delta}(y)|y|^{\rho}dy.
\end{align}
From (\ref{delta}) and (\ref{alpha}), it follows that
\begin{align*}
\int_{\mathbb{R}}(p_{t+\delta}(x-y)-p_{t}(x-y))w(y)dy\leq 2\int_{0}^{\infty}\frac{\delta y^{\rho}}{(\delta^{1/\alpha}+y)^{1+\alpha}}dy=2\delta^{\rho/\alpha}\int_{0}^{\infty}\frac{y^{\rho}}{(1+y)^{1+\alpha}}dy.
\end{align*}
\end{proof}

\subsection{Difference in the spatial variable}
For a random variable $X\in L^{k}(P)$,  define $\|X\|_{L^{k}(P)}=(E(|X|^{k}))^{1/k}$. For simplicity, we write $\|\cdot\|_{k}$ instead of $\|\cdot\|_{L^{k}(P)}$.
We will estimate the spatial and the time difference of the following stochastic integral $I$ in this and next subsection.
For $t\in [0,T]$, $x,y,z\in \mathbb{R}$,
Define
\begin{align}
I_{\beta,t}(x)=\int_{0}^{t}\int_{\mathbb{R}}p_{t-s}(z-x))\sigma(u_{\beta,s}(z))\eta_{\beta}(ds,dz),
\end{align}
and denote
\begin{align*}
A_{s}(x,y)&=\sigma(u_{\beta,s}(x))\sigma(u_{\beta,s}(y)),\\
B_{s}(r)&=p_{t-s}(r-x)-p_{t-s}(r-y).
\end{align*}
For all $k\geq 2$, the difference in the spatial variable is
\begin{align*}
E\big(|I_{\beta,t}(x)-I_{\beta,t}(y)|^{k}\big)=E\big(\big|\int_{0}^{t}\int_{\mathbb{R}}(p_{t-s}(z-x)-p_{t-s}(z-y))\sigma(u_{\beta,s}(z))\eta_{\beta}(ds,dz)\big|^{k}\big).
\end{align*}

\begin{theorem}\label{spatial varibale}
For all $t\in [0,T], \ x,y\in \mathbb{R}$,
\begin{align}
E\big(|I_{\beta,t}(x)-I_{\beta,t}(y)|^{k}\big)\leq C|x-y|^{\frac{\alpha bk}{2}},
\end{align}
where $C$ is a constant, $b\in (0, 1-\frac{1}{\alpha})$, $\alpha\in (1,2]$, $\beta\in (0,1)$.
\end{theorem}

\begin{proof}
We apply the Cauchy-Schwarz inequality to bound $E(|\sigma(u_{\beta,s}(x))\sigma(u_{\beta,s}(y))|^{k/2})$ by $\sup_{x \in\mathbb{R}}E|\sigma(u_{\beta,s}(x))|^{k}$. Similar to the proof of Theorem 13 of \cite{Robert}, we can show that $\sup_{s\in [0,T]}\sup_{x \in\mathbb{R}}E|u_{\beta,s}(x)|^{k}< \infty$. Then by (A2), we obtain $\sup_{s\in [0,T]}\sup_{x \in\mathbb{R}}E|\sigma(u_{\beta,s}(x))|^{k}<\infty$.
It is easy to show that $p_{r}\ast f$ is positive definite and continuous for all $r>0$, then $\sup_{z\in \mathbb{R}}(p_{t-s}\ast f_{\beta})(z)=(p_{t-s}\ast f)(0)$. Since $(B_{s}\ast f_{\beta})(w)\leq 2\sup_{z\in \mathbb{R}}(p_{t-s}\ast f)(z)$, we get $\sup_{z\in \mathbb{R}}(p_{t-s}\ast f_{\beta})(z)\leq 2(p_{t-s}\ast f)(0)$. By Burkholder inequality (\cite{DKMNX}, Theorem 5.27), Minkowski integral inequality (\cite{Stein}, Appendice A.1) and Lemma \ref{space} we have
\begin{align}\label{spatial estimate}
&E\big(|I_{\beta,t}(x)-I_{\beta,t}(y)|^{k}\big)\nonumber\\
&\leq c_{k}E\big(\big|\int_{0}^{t}\int_{\mathbb{R}}\int_{\mathbb{R}}f_{\beta}(z-w)B_{s}(z)B_{s}(w)A_{s}(z,w)dsdzdw\big|^{k/2}\big)\nonumber\\
&\leq c_{k}\big|\int_{0}^{t}\sup_{x\in \mathbb{R}}\|\sigma(u_{\beta,s}(x))\|_{k}^{2}\int_{\mathbb{R}}\int_{\mathbb{R}}f_{\beta}(z-w)|B_{s}(z)||B_{s}(w)|dsdzdw\big|^{k/2}\nonumber\\
&\leq C\big|\int_{0}^{t}\int_{\mathbb{R}}\int_{\mathbb{R}}f_{\beta}(z-w)|B_{s}(z)||B_{s}(w)|dsdzdw\big|^{k/2}\nonumber\\
&\leq C\big|\int_{0}^{t}(p_{t-s}\ast f_{\beta})(0)ds\int_{\mathbb{R}}|p_{t-s}(z-x)-p_{t-s}(y-x)|dz\big|^{k/2}\nonumber\\
&\leq C\big|\int_{0}^{t}(p_{t-s}\ast f_{\beta})(0)\big(\frac{|x-y|}{(t-s)^{1/\alpha}}\wedge 1\big)ds\big|^{k/2}.
\end{align}
Since $r\wedge 1\leq r^{\alpha b}$ for all $r>0$ and $b\in (0,\frac{1}{\alpha})$, by (\ref{spatial estimate}), we have
\begin{align}\label{constant}
E\big(|I_{\beta,t}(x)-I_{\beta,t}(y)|^{k}\big)\leq C|x-y|^{\frac{\alpha bk}{2}}\big|\int_{0}^{t}(p_{t-s}\ast f_{\beta})(0)(t-s)^{-b}ds\big|^{k/2}.
\end{align}
By (\ref{alpha}), we have
\begin{align}\label{sup}
(p_{t-s}\ast f_{\beta})(0)&=c_{1-\beta}\int_{\mathbb{R}}\frac{1}{|x|^{\beta}}p_{t-s}(x)dx\leq c_{1-\beta}\int_{\mathbb{R}}\frac{1}{|x|^{\beta}}\cdot \frac{t-s}{((t-s)^{1/\alpha}+|x|)^{1+\alpha}}dx\nonumber\\
&\leq c_{1-\beta}(t-s)^{-\beta/\alpha}\int_{\mathbb{R}}\frac{1}{|r|^{\beta}(1+|r|)^{1+\alpha}}dr\leq C(t-s)^{-\beta/\alpha}.
\end{align}
Put (\ref{sup}) into (\ref{constant}) to get
\begin{align}\label{spatial bound}
E\big(|I_{\beta,t}(x)-I_{\beta,t}(y)|^{k}\big)\leq C|x-y|^{\frac{\alpha bk}{2}}\big|\int_{0}^{t}(t-s)^{-\frac{\beta}{\alpha}-b}ds\big|^{k/2}.
\end{align}
Since $\beta\in (0,1)$ and $\alpha\in (1,2]$, we can choose $b\in (0, 1-\frac{1}{\alpha})\subset (0, \frac{1}{\alpha})$ to guarantee that $\big|\int_{0}^{t}(t-s)^{-\frac{\beta}{\alpha}-b}ds\big|<\infty$.
Therefore we obtain
\begin{align}\label{spatial continuity}
E\big(|I_{\beta,t}(x)-I_{\beta,t}(y)|^{k}\big)\leq C|x-y|^{\frac{\alpha bk}{2}}, \ t\in [0,T], \ x,y\in \mathbb{R}.
\end{align}
\end{proof}

\subsection{Difference in the time variable}
For all $k\geq 2$, the difference in the time variable is
\begin{align*}
&E(|I_{\beta,t+\delta}-I_{\beta,t}|^{k})\\
&=E(|\int_{0}^{t+\delta}\int_{\mathbb{R}}p_{t+\delta-s}(z-x)\sigma(u_{\beta,s}(z))\eta_{\beta}(ds,dz)-\int_{0}^{t}\int_{\mathbb{R}}p_{t-s}(z-x)\sigma(u_{\beta,s}(z))\eta_{\beta}(ds,dz)|^{k})
\end{align*}

\begin{theorem}\label{time variable}
For all $\delta>0$, $t\in [0,T]$, $x, y\in \mathbb{R}$,
\begin{align}
E(|I_{\beta,t+\delta}(x)-I_{\beta,t}(x)|^{k})\leq C\delta^{\frac{(\alpha-\beta)k}{2\alpha}},
\end{align}
where $C$ is a constant, $\alpha\in (1,2]$, $\beta\in (0, \frac{\alpha}{2}]$.
\end{theorem}

\begin{proof}
By the elementary inequality $|a+b|^{k}\leq 2^{k}|a|^{k}+2^{k}|b|^{k}$, we have
\begin{align}
&E\big(|I_{\beta,t+\delta}(x)-I_{\beta,t}(x)|^{k})\nonumber\\
&\leq 2^{k}E\big(|\int_{0}^{t+\delta}\int_{\mathbb{R}}(p_{t+\delta-s}(z-x)-p_{t-s}(z-x))\sigma(u_{\beta,s}(z))\eta_{\beta}(ds,dz)|^{k}\big)\nonumber\\
&\quad+2^{k}E\big(|\int_{t}^{t+\delta}\int_{\mathbb{R}}(p_{t+\delta-s}(z-x)\sigma(u_{\beta,s}(z))\eta_{\beta}(ds,dz)|^{k}\big)\nonumber\\
&=I_{1}+I_{2}.
\end{align}
For $I_{2}$, by the same technique as in the proof of Theorem \ref{spatial varibale}, we have
\begin{align}\label{second}
I_{2}&\leq C\big(\int_{t}^{t+\delta}\sup_{x\in \mathbb{R}}\|\sigma(u_{\beta,s}(x))\|_{k}^{2}\int_{\mathbb{R}}\int_{\mathbb{R}}f_{\beta}(z-w)p_{t+\delta-s}(z-x)p_{t+\delta-s}(w-x)dsdzdw\big)^{k/2}\nonumber\\
&\leq C\big(\int_{t}^{t+\delta}\int_{\mathbb{R}}\int_{\mathbb{R}}f_{\beta}(z-w)p_{t+\delta-s}(z-x)p_{t+\delta-s}(w-x)dsdzdw\big)^{k/2}.
\end{align}
Denote by $\mathcal{S}(\mathbb{R})$ the Schwartz space of rapid decreasing test-functions from $\mathbb{R}$ to $\mathbb{R}$, by elementary properties of convolution and Fourier transform, the following holds (see formula (10) in \cite{Robert}):
\begin{align*}
\int_{\mathbb{R}}\int_{\mathbb{R}}\varphi(x)f_{\beta}(x-y)\varphi(y)dxdy=\int_{\mathbb{R}}f_{\beta}(x)(\varphi\ast \tilde{\psi})(x)dx=\int_{\mathbb{R}}g_{1-\beta}(\xi)|\mathcal{F}\varphi(\xi)|^{2}d\xi,
\end{align*}
for all $\varphi, \psi\in \mathcal{S}(\mathbb{R})$, where $\tilde{\psi}$ is defined by $\tilde{\psi}(x)=\psi(-x)$.
Then by change of variables,
\begin{align}\label{fractional transform}
\int_{\mathbb{R}}\int_{\mathbb{R}}p_{t+\delta-s}(z-x)f_{\beta}(z-w)p_{t+\delta-s}(w-x)dzdw=\int_{\mathbb{R}}g_{1-\beta}(\xi)|\mathcal{F}p_{t+\delta-s}(\xi)|^{2}d\xi.
\end{align}
Recall the Fourier transform of fractional heat kernel $p_{t}(x)$ (see formula (3) in \cite{KT}) and note (\ref{Fourier transform}), we have
\begin{align}\label{fthk}
\mathcal{F}p_{t-s}(\xi)=e^{-(t-s)(2\pi |\xi|)^{\alpha}}.
\end{align}
Thus,
\begin{align}\label{fractional convoution}
\int_{\mathbb{R}}g_{1-\beta}(\xi)|\mathcal{F}p_{t-s}(\xi)|^{2}d\xi&=\int_{\mathbb{R}}g_{1-\beta}(\xi)e^{-(t-s)2^{\alpha+1}\pi^{\alpha}|\xi|^{\alpha}}d\xi\nonumber\\
&=\int_{\mathbb{R}}|\xi|^{\beta-1}e^{-(t-s)2^{\alpha+1}\pi^{\alpha}|\xi|^{\alpha}}d\xi\nonumber\\
&=2\int_{0}^{\infty}\xi^{\beta-1}e^{-(t-s)2^{\alpha+1}\pi^{\alpha}\xi^{\alpha}}d\xi\nonumber\\
&=\frac{2}{\alpha}\int_{0}^{\infty}e^{-(t-s)2^{\alpha+1}\pi^{\alpha}r}r^{\frac{\beta}{\alpha}-1}dr\nonumber\\
&=\frac{2}{\alpha(2^{\alpha+1}\pi^{\alpha}(t-s))^{\beta/\alpha}}\int_{0}^{\infty}e^{-z}z^{\frac{\beta}{\alpha}-1}dz\nonumber\\
&=\frac{2\Gamma(\frac{\beta}{\alpha})}{\alpha(2^{\alpha+1}\pi^{\alpha}(t-s))^{\beta/\alpha}}.
\end{align}
By (\ref{second}), (\ref{fractional transform}) and (\ref{fractional convoution}), we get
\begin{align}\label{first bound}
I_{2}\leq C\big(\int_{t}^{t+\delta}(t+\delta-s)^{-\beta/\alpha}ds\big)^{k/2}=C\delta^{\frac{(\alpha-\beta)k}{2\alpha}}.
\end{align}
For $I_{1}$, by the similar argument and note that Lemma \ref{time} and (\ref{sup}), we have
\begin{align}\label{first}
I_{1}&\leq C\big(\int_{0}^{t}\sup_{x\in \mathbb{R}}\|\sigma(u_{\beta,s}(x))\|_{k}^{2}\ (f_{\beta}\ast p_{t-s})(0)\int_{\mathbb{R}}p_{t+\delta-s}(z)-p_{t-s}(z)dzds\big)^{k/2}\nonumber\\
&\leq C\big(\int_{0}^{t}s^{-\beta/\alpha}(\log(s+\delta)-\log(s))ds\big)^{k/2}.
\end{align}
By integrating by parts,
\begin{align}\label{integral by parts}
&\int_{0}^{t}s^{-\beta/\alpha}(\log(s+\delta)-\log(s))ds\nonumber\\
&=\frac{\alpha}{\alpha-\beta}\log(1+\delta/t)t^{(\alpha-\beta)/\alpha}+\frac{\alpha}{\alpha-\beta}\int_{0}^{t}s^{(\alpha-\beta)/\alpha}\frac{\delta}{s(s+\delta)}ds\nonumber\\
&=I_{3}+I_{4}.
\end{align}
For $I_{4}$,
\begin{align}\label{delta bound}
\frac{\alpha}{\alpha-\beta}\int_{0}^{t}s^{(\alpha-\beta)/\alpha}\frac{\delta}{s(s+\delta)}ds&=\big(\frac{\alpha}{\alpha-\beta}\big)^{2}\int_{0}^{\frac{(\alpha-\beta)t}{\alpha}}\frac{\delta}{\mu^{\alpha/(\alpha-\beta)}+\delta}d\mu\nonumber\\
&=\big(\frac{\alpha}{\alpha-\beta}\big)^{2}\delta^{\frac{(\alpha-\beta)}{\alpha}}\int_{0}^{\frac{(\alpha-\beta)t}{\alpha}\delta^{\frac{(\beta-\alpha)}{\alpha}}}\frac{1}{1+\nu^{\alpha/(\alpha-\beta)}}d\mu\nonumber\\
&\leq \big(\frac{\alpha}{\alpha-\beta}\big)^{2}\delta^{\frac{(\alpha-\beta)}{\alpha}}\int_{0}^{\infty}\frac{1}{1+\nu^{\alpha/(\alpha-\beta)}}d\mu \nonumber\\
&\leq C\delta^{\frac{(\alpha-\beta)}{\alpha}}.
\end{align}

Next, we will prove that for any $\mu>0$, $r\in [\frac{1}{2},1]$,
\begin{align}\label{log bound}
0<\log(1+\mu)\leq \mu^{r}.
\end{align}
For $r=1$, it is obvious that $0<\log(1+\mu)\leq \mu$. In the following, we only consider the case $r\in [\frac{1}{2},1)$. For $\mu>0$, let $h(\mu)=\log(1+\mu)-\mu^{r}$. Then
\begin{align}\label{derivative}
h'(\mu)=\frac{\mu^{1-r}-(1+\mu)r}{(1+\mu)\mu^{1-r}}.
\end{align}
Let $l(\mu)=\mu^{1-r}-(1+\mu)r$. It is easy to get the maximum $l_{max}=l(\mu)|_{\mu=(\frac{1-r}{r})^{\frac{1}{r}}}=\frac{r^{2}}{1-r}\big[(\frac{1-r}{r})^{\frac{1}{r}}-\frac{1-r}{r}\big]$. Since $r\in [\frac{1}{2},1)$, we get  $l_{max}\leq 0$. By (\ref{derivative}), $h'(\mu)\leq 0$. This together with $h(0)=0$ yield that $h(\mu)\leq 0$. Therefore (\ref{log bound}) holds.
Since $\beta\in (0, \frac{\alpha}{2}]$,  then $\frac{\alpha-\beta}{\alpha}\in [\frac{1}{2},1)$.  By (\ref{log bound}),  we have
\begin{align}\label{logdelta}
0<\log(1+\delta/t)\leq (\delta/t)^{(\alpha-\beta)/\alpha}.
\end{align}
Thus, for $I_{3}$, we get
\begin{align}\label{log}
I_{3}\leq C\delta^{\frac{(\alpha-\beta)}{\alpha}}.
\end{align}
By (\ref{first}), (\ref{integral by parts}), (\ref{delta bound}) and (\ref{log}),  we have
\begin{align}
E(|I_{\beta,t+\delta}(x)-I_{\beta,t}(x)|^{k})\leq C\delta^{\frac{(\alpha-\beta)k}{2\alpha}}.
\end{align}

\end{proof}
\subsection{Spatial and temporal H$\ddot{\mbox{o}}$lder continuity}

\begin{theorem}
For all $k\geq 2$, $\alpha\in (1,2]$, $\beta\in (0,\frac{\alpha}{2}]$, $\rho\in (0,1)$, $b\in (0,1-\frac{1}{\alpha})$ and $x,y\in \mathbb{R}$,
\begin{align}
E\big(|u_{\beta,t}(x)-u_{\beta,s}(y)|\big)^{k}\leq C\big(|x-y|^{kc}+|t-s|^{kd}\big),
\end{align}
where $C_{1}$, $C_{2}$ are positive constants, $c\in (0, \frac{\alpha b}{2}\wedge \rho)$, $d\in (0,\frac{(\alpha-\beta)}{2\alpha}\wedge \frac{\rho}{\alpha})$.
\end{theorem}
\begin{proof}
By Lemma \ref{holder}, Theorem \ref{spatial varibale} and Theorem \ref{time variable}, we can draw the conclusion.
\end{proof}

\section{Acknowledgements}

This work is partial supported by National Natural Science Foundation of China under the contract
No.11571269, China Postdoctoral Science Foundation Funded Project under contracts No.2015M572539 and No.2016T90899 and Shaanxi Province Postdoctoral Science Foundation Funded Project.\


\end{document}